\documentclass{article}
\usepackage{fullpage}
\usepackage{graphicx} 
\usepackage{amsmath}
\usepackage{bbm}
\usepackage{amsfonts}
\usepackage{amssymb}
\usepackage{epsfig}
\usepackage{amsthm}
\usepackage{url}
\usepackage{cite}

\begin{document}
\title{A Two-Step High-Order Compact Scheme for the Laplacian Operator and its Implementation in an Explicit Method for Integrating the Nonlinear Schr{\"o}dinger Equation.}
\author{R. M. Caplan\footnote{Corresponding author.  Present address: Predictive Science Inc. 9990 Mesa Rim Rd, Suite 170, San Diego, CA 92121. email: caplanr@predsci.com, phone: 858-225-2314, \texttt{URL}: http://www.predsci.com} and R. Carretero-Gonz{\'a}lez\\[1.0ex]
Nonlinear Dynamical Systems Group\footnote{\texttt{URL}: http://nlds.sdsu.edu},
Computational Science Research Center\footnote{\texttt{URL}: http://www.csrc.sdsu.edu}, and\\
Department of Mathematics and Statistics,
San Diego State University,\\
San Diego, California 92182-7720, USA\\
}
\date{\today}
\maketitle
\begin{abstract}
We describe and test an easy-to-implement two-step high-order compact (2SHOC) scheme for the Laplacian operator and its implementation into an explicit finite-difference scheme for simulating the nonlinear Schr{\"o}dinger equation (NLSE).  Our method relies on a compact `double-differencing' which is shown to be computationally equivalent to standard fourth-order non-compact schemes.  Through numerical simulations of the NLSE using fourth-order Runge-Kutta, we confirm that our scheme shows the desired fourth-order accuracy.  A computation and storage requirement comparison is made between the 2SHOC scheme and the non-compact equivalent scheme for both the Laplacian operator alone, as well as when implemented in the NLSE simulations.  Stability bounds are also shown in order to get maximum efficiency out of the method.  We conclude that the modest increase in storage and computation of the 2SHOC schemes are well worth the advantages of having the schemes compact, and their ease of implementation makes their use very useful for practical implementations.
\\
\\
\textit{Key Words:} High-order Compact Scheme, Nonlinear Schr{\"o}dinger equation.
\end{abstract}

\section{Introduction}
\label{s:intro}
The nonlinear Schr{\"o}dinger equation (NLSE) is a universal model describing the evolution and propagation of complex field envelopes in nonlinear dispersive media.  As such, it is used to describe many physical systems including the evolution of water waves, nonlinear optics, thermodynamic pulses, nonlinear waves in fluid dynamics, and waves in semiconductors \cite{NLSE_nlpdebook, NLSE_selffocus,NLSE_OPTVA}.  The general form of the NLSE can be written as
\begin{equation}
\label{NLSE}
i\frac{\partial \Psi}{\partial t} + a\nabla^2\Psi - V({\bf r})\Psi + s|\Psi|^2 \Psi = 0,
\end{equation}
where $\Psi({\bf r},t) \in \mathbbm{C} $ is the value of the wavefunction, $\nabla^2$ is the Laplacian operator, and where $a>0$ and $s$ are parameters defined by the system being modeled.  $V(\bf{r})$ is an external potential term, which when included, makes Eq.~(\ref{NLSE}) known as the Gross-Pitaevskii equation \cite{BEC_RCbook}.

The nonlinearity in the NLSE allows the prediction and description 
of important, and experimentally relevant, nonlinear effects and 
nonlinear waves, such as solitons and vortices~\cite{BEC_RCbook}.
In the context of nonlinear optics, this nonlinearity emerges from
the nonlinear response of the optical medium. The cubic-type nonlinearity
in this context is the so-called Kerr nonlinearity of the medium 
\cite{NLKerr,NLScrys}.
Another context where the NLSE serves as a prototypical model is in the
realm of Bose-Einstein condensates (BECs)~\cite{BEC_RCbook}.
In this context, the nonlinearity takes into account the {\em mean-field}
interaction of the atoms in the condensate.
The sign of the nonlinearity in the NLSE (i.e., the sign of $s$)
defines the type of interactions. For example, in the nonlinear
optics application a self-focusing (self-defocusing) medium corresponds
to $s>0$ ($s<0$); while in the BEC context, attractive (repulsive)
interatomic interactions correspond to $s>0$ ($s<0$).
In the case of self-focusing or attracting nonlinearity the 
prototypical (nonlinear) solutions of the NLSE are in the form
of {\em bright} solitons in one dimension; while in two or three
dimensions bright vortices and vortex rings do exists but are typically
unstable.
On the other hand, the case of self-defocusing or repulsive
nonlinearity bears {\em dark} (nonlinear) solutions that in
one dimension take the shape of dark (non-moving) or grey
(moving) solitons which in two and three dimensions correspond 
to, typically stable, (dark) vortices and rings, respectively.

For almost every nonlinear partial differential equation (PDE) (including the NLSE), the interesting solutions require simulations using numerical methods.   A standard methodology includes discretizing the solution in space by using finite-difference approximation schemes for the required derivatives.  A \emph{high-order} scheme is one which exhibits greater than second-order ($O(h^2)$, where $h$ is the grid spacing) accuracy.  While there are a large number of scenarios where higher-order schemes are a necessity due to the desired accuracy of the simulations, often the higher accuracy of high-order schemes is unnecessary, and second-order accuracy is sufficient for the problem.  However, in such a case, using high-order schemes is still desirable because they allow one to simulate a solution with many fewer grid points, while maintaining the same accuracy as a second-order scheme.  In those situations, the desired grid point size is based on the ability to resolve the structure of the solution, and not on the accuracy of computation.

Higher-order finite-difference schemes are typically achieved by computing derivatives with a wider scheme stencil.  This adds a difficulty near the boundary, as one must be able to calculate the inner point near the boundary with the same accuracy as the internal scheme which can be complicated to implement.  If this is not done, the entire simulation can eventually become lower-order. Another disadvantage of wide-stencils is that they cause many parallel implementations to be more difficult to realize, as different compute nodes must share or transfer more boundary values to their neighboring nodes.  Besides the added complexity of the codes, the extra communication/global-memory-access (in the case of graphical processing units) reduces the overall parallel performance \cite{HOC_PARPER}.  

For the aforementioned reasons, there is great interest in high-order compact (HOC) schemes.  These are finite-difference schemes which exhibit higher-order accuracy but still only rely on the closest neighboring points for computations. HOC schemes have been developed for various multidimensional steady-state PDEs including the convection-diffusion equation \cite{HOC_3D_SS_ConvDiff}, Poisson's equation \cite{HOC_3DPOSS}, the stream-function vorticity form of the Navier-Stokes equations \cite{HOC_2D_vel}, as well as generalized linear elliptical PDEs with variable coefficients \cite{HOC_3D_MAPLE} to name a few.  Many HOC implementations of time-dependent PDEs have also been formulated including Burger's equation \cite{HOC_Burgers}, the wave equation \cite{HOC_TWAVE}, the Euler equations \cite{HOC_Euler_HOC} and the time-dependent convection-diffusion equation \cite{HOC_TIME,HOC_2DConDiffT} as well as others.  Time-dependent HOC schemes have also been developed for the one- and two-dimensional linear \cite{HOC_LSE2D_BVM,HOC_LSE2D_OLD} and one-dimensional nonlinear Schr{\"o}dinger equations \cite{HOC_1DNLS,HOC_1DNLSSS}.  One drawback of HOC schemes are that their formulations are typically specific for the model equation being used, and requires re-deriving the schemes for each different PDE to be simulated (see Ref.~\cite{HOC_3D_MAPLE}, where the authors developed a Maple program to symbolically generate a HOC scheme for three-dimensional linear elliptical models).  

The time-dependent HOC schemes developed so far are all implicit (even if the time-stepping is done with an otherwise explicit scheme), requiring solving a linear system in each time-step, as well as iterative processes when simulating nonlinear PDEs such as the NLSE (see however Ref.~\cite{HOC_1DNLS} where the author's second implicit scheme implementation was able to be formulated to get around this requirement).  The computational and storage requirements for implicit schemes can be prohibitive in large multidimensional settings.   They are also, in general, difficult to optimally parallelize, especially in higher dimensions.  

For these reasons, we wish to develop HOC formulations that are {\em fully explicit} time-dependent schemes (in our case, for the NLSE).  To do this, we formulate a two-step procedure in computing the spatial derivatives, where each individual step is a compact computation.  The first step of the 2SHOC scheme computes the standard second-order finite-difference approximation to the derivatives of the Laplacian, while the second step uses these computed derivative values to compactly approximate the Laplacian to fourth-order accuracy.  Since the 2SHOC scheme computes the Laplacian independent of the governing PDE, it can be used in a multitude of time-dependent and time-independent PDEs which contain the Laplacian operator, increasing its generality.  An explicit 2SHOC scheme for the NLSE can be formed by utilizing standard explicit ordinary differential equation (ODE) solvers.  For the scheme described in this paper, we use the classic fourth-order Runge-Kutta (RK4)\cite{RK4}. 

The paper is organized as follows.  In Section~\ref{s:form} we show the formulation of the 2SHOC scheme for the Laplacian operator in one, two and three dimensions.  Then, in Section~\ref{s:2shocrk4}, we use the scheme to form an explicit method for simulating the NLSE.  In Section~\ref{s:num}, we show numerical tests of each scheme from Section~\ref{s:2shocrk4} and confirm their accuracy.  In Section~\ref{s:ccompare} we compare the storage and computation requirements for the 2SHOC schemes versus the standard fourth-order non-compact equivalent schemes.  We conclude in Section~\ref{s:con}.

\section{Formulation of 2SHOC Schemes for the Laplacian Operator}
\label{s:form}
It is well known that one can derive a fourth-order accurate finite-difference scheme for the second spatial derivative of a function by applying a `double-differencing' approach.  This works by noting that the first truncation term in the standard central-difference scheme contains the fourth spatial derivative.  A second-order approximation to the fourth derivative can be obtained by applying a central-difference operator to the result of applying a central-difference operator to the function.  When multiplied by the $h^2$ (where $h$ is the spatial step-size of the computational grid) in the truncation term, the error in the truncation term becomes $O(h^4)$, and therefore, the resulting scheme becomes a fourth-order accurate approximation to the second derivative.  When these two steps are algebraically combined and simplified, one yields the standard non-compact fourth-order stencil.  Due to the increase in the number of computations and storage, the `double-differencing' procedure is not actually implemented, rather, the resulting non-compact stencil is used directly. 

However, when the advantages of using a compact scheme are considered, the numerical implementation of the  `double-differencing' procedure becomes an overlooked, viable alternative to other more complicated HOC schemes.  In time-dependent systems, the sequential nature of the `double-differencing' procedure makes it best suited for use with explicit time-stepping schemes.  This may be a reason why this simple compact solution has been overlooked, as many high-order numerical schemes for time-dependent PDEs are typically implicit.

\subsection{One-Dimensional Formulation}
\label{sec:2shoc1d}
Although the methodology for the 2SHOC scheme is to use the `double-differencing' approach, it is nevertheless worthwhile to show that the scheme can also be formed as a result of attempting to make other time-dependent HOC scheme formulations explicit.  In Ref.~\cite{HOC_TIME}, the authors formulate an implicit time-dependent HOC scheme for the transient convection-diffusion equation based on their method from Ref.~\cite{HOC_3DPOSS}.  We can obtain a compact high-order scheme for the time-dependent NLSE following the same methodology, and with a modification, can transform it into an explicit scheme which is equivalent to using the 2SHOC scheme for the Laplacian operator in the NLSE.

In one dimension, we discretize the wavefunction of the NLSE as $\Psi(x_i,t_n) \equiv \Psi_i^n$ with a grid spacing of size $h$ ($x_i = x_{\mbox{\scriptsize min}} + hi$) and a time-step of size $k$ ($t_n = kn$).   The Laplacian operator applied to $\Psi$ in one dimension can be represented as
\begin{equation}
\label{dx2full}
\nabla^2 \Psi_i = \left. \frac{\partial^2 \Psi}{\partial x^2}\right|_i = \delta_x^2 \Psi_i - \frac{h^2}{12}\left(\left. \frac{\partial^4 \Psi}{\partial x^4}\right|_i\right) + O(h^4),
\end{equation}
where the central-difference operator $\delta_x^2$ is defined as
\begin{equation}
\label{dx2}
\delta_x^2 \,\Psi_i = \frac{1}{h^2}(\Psi_{i+1} - 2\Psi_i + \Psi_{i-1}).
\end{equation}

As per the technique of Ref.~\cite{HOC_TIME}, to formulate a high-order compact scheme, the NLSE of Eq.~(\ref{NLSE}) is differentiated in space twice in order to form an expression for the fourth derivative in Eq.~(\ref{dx2full}) which yields
\begin{equation}
\label{dx4}
\left. \frac{\partial^4 \Psi}{\partial x^4}\right|_i = -\frac{1}{a}\; \delta_x^2\left[i\left. \frac{\partial \Psi}{\partial t}\right|_i + \left(s|\Psi_i|^2 - V(x_i)\right)\Psi_i\right] + O(h^2).
\end{equation}
Due to the $h^2$ in the truncation term of Eq.~(\ref{dx2full}), when Eq.~(\ref{dx4}) is inserted into  Eq.~(\ref{dx2full}), the resulting approximation of the Laplacian becomes $O(h^4)$.  Inserting Eqs.~(\ref{dx2full}) and (\ref{dx4}) into Eq.~(\ref{NLSE}) yields the fourth-order in space semi-discrete HOC scheme
\begin{alignat}{2}
\label{nlsimphoc1d}
\left. \frac{\partial \Psi}{\partial t}\right|_i = &i [ a\left(\delta_x^2\,\Psi_i + \frac{h^2}{a\,12}\delta_x^2\left[i\,\left. \frac{\partial \Psi}{\partial t}\right|_i + \left(s\,|\Psi_i|^2 - V(x_i)\right)\Psi_i\right]_i \right)
\\
&+ \left(s\,|\Psi_i|^2-V(x_i)\right)\Psi_i]. \notag
\end{alignat}
Due to the central-difference operator operating on the temporal derivative $\partial \Psi/\partial t$, the resulting scheme of Eq.~(\ref{nlsimphoc1d}) will be implicit even if the time-stepping is chosen to be otherwise explicit (for example, forward differencing).

In order to retain a fully explicit scheme, a two-step approach is used.  First, an approximation of $\partial \Psi_i/\partial t$ is made and then is inserted into Eq.~(\ref{nlsimphoc1d}).  The approximation of $\partial \Psi_i/\partial t$  must be $O(h^2)$ for the scheme to retain its fourth-order accuracy.  The most straight-forward way to make such an approximation is to apply a standard second-order central differencing to the NLSE yielding
\begin{equation}
\label{uttmp}
T_i = \left. \frac{\partial \Psi}{\partial t}\right|_i + O(h^2) = i\left[a\delta_x^2\,\Psi_i + \left(s\,|\Psi_i|^2-V(x_i)\right)\Psi_i\right].
\end{equation}
Once computed, $T_i$ is inserted into Eq.~(\ref{nlsimphoc1d}), and the central difference operator can then be applied to $T_i$ as
\[
\delta_x^2 \,T_i = \frac{T_{i+1} - 2\,T_i + T_{i-1}}{h^2}.
\]
It is important that the boundary conditions for Eq.~(\ref{uttmp}) be computed to at least $O(h^2)$ accuracy since the boundary points will be used for interior calculations of the spatial derivative of Eq.~(\ref{nlsimphoc1d}) (see Sec.~\ref{s:bc} for details on boundary conditions for the 2SHOC schemes).  Inserting Eq.~(\ref{uttmp}) into Eq.~(\ref{nlsimphoc1d}) yields the semi-discrete equation
It is important that the boundary conditions for Eq.~(\ref{uttmp}) be computed to at least $O(h^2)$ accuracy since the boundary points will be used for interior calculations of the spatial derivative of Eq.~(\ref{nlsimphoc1d}) (see Sec.~\ref{s:bc} for details on boundary conditions for the 2SHOC schemes).  Inserting Eq.~(\ref{uttmp}) into Eq.~(\ref{nlsimphoc1d}) yields the semi-discrete equation
\begin{alignat}{2}
\label{ehoc}
\left. \frac{\partial \Psi}{\partial t}\right|_i = &i [a\left(\delta_x^2\,\Psi_i + \frac{h^2}{a\,12}\delta_x^2\left[i\,T_i + \left(s\,|\Psi_i|^2 - V(x_i)\right)\Psi_i\right]_i \right)
\\
& + \left(s\,|\Psi_i|^2-V(x_i)\right)\Psi_i] + O(h^4). \notag
\end{alignat}
Any desired explicit ODE solver (that is stable for the problem with its parameters) can then be used to integrate Eq.~(\ref{ehoc}).

Algebraically combining the two stages of Eq.~(\ref{uttmp}) and Eq.~(\ref{ehoc}) together yields 
\begin{equation}
\label{ddd}
\left. \frac{\partial \Psi}{\partial t}\right|_i = i\left[a\left(\delta_x^2\,\Psi_i - \frac{h^2}{12}\delta_x^2 (\delta_x^2 \Psi_i)\right)+ \left(s\,|\Psi_i|^2-V(x_i)\right)\Psi_i\right] + O(h^4),
\end{equation}
and therefore the two-stage HOC scheme described in Eq.~(\ref{uttmp}) and Eq.~(\ref{ehoc}) is computationally equivalent to simply taking the central-difference of the central-difference to approximate the fourth derivative truncation term in the Laplacian.  Therefore, we have recovered the 2SHOC scheme approach for approximating the Laplacian which does not depend on the other terms of the governing equation.  Therefore the 2SHOC is a stand-alone scheme for the Laplacian operator which can be used in multiple governing equations. 

After collecting terms and simplifying the 2SHOC `double-differencing', the resulting two-step scheme for the one-dimensional Laplacian is given by
\begin{alignat}{3}
&1) \qquad &D_i &= \frac{1}{h^2}\left(\Psi_{i+1} - 2\Psi_i + \Psi_{i-1}\right), \label{2shoc1d}\\
&2) \qquad &\nabla^2\Psi_i &\approx \frac{7}{6}D_i - \frac{1}{12}\left(D_{i+1} + D_{i-1}\right).\label{2shoc1d2}
\end{alignat}

When the two steps of Eqs.~(\ref{2shoc1d}) and (\ref{2shoc1d2}) are combined algebraically and simplified, the standard five-point non-compact fourth-order finite-difference approximation is recovered:
\[
\nabla^2\Psi_i = \left. \frac{\partial^2 \Psi}{\partial x^2} \right|_i = \frac{-\Psi_{i+2} +16\Psi_{i+1} -30\Psi_i +16\Psi_{i-1} - \Psi_{i-2}}{12h^2} + O(h^4).
\]

A potential drawback in using the 2SHOC scheme is that it requires extra storage space (the $D$ array) and more computations than using the standard fourth-order five-point stencil.  However, as will be discussed in Sec.~\ref{s:ccompare}, the compact scheme's advantages can outweigh this deficiency. 

\subsection{Two-Dimensional Formulation}\label{sec:2shoc2d}
In two dimensions, the Laplacian operator applied to $\Psi$ at grid location $(i,j)$ can be represented as
\begin{equation}
\label{dx2full2d}
\begin{tabular}{c} 
$\nabla^2 \Psi_{i,j} = \dfrac{1}{h^2}$
\begin{tabular}{|c|c|c|} \hline
  &  1 &   \\ \hline
 1 &-4 &  1 \\ \hline
  &  1 &  \\ \hline
\end{tabular}
$\Psi_{i,j} - \dfrac{h^2}{12}\left( \left. \dfrac{\partial^4 \Psi}{\partial x^4}\right|_{i,j} + \left. \dfrac{\partial^4 \Psi}{\partial y^4}\right|_{i,j} \right) + O(h^4).$
\end{tabular}
\end{equation}
Unlike the one-dimensional case, there are additional compact grid points which are not being used in Eq.~(\ref{dx2full2d}) (the four corner points).  It is known that these points can be added to make a more accurate nine-point Laplacian operator given as \cite{HOC_ninepointBOOK}
\begin{equation}
\label{9ptlap}
\begin{tabular}{c} 
$\nabla^2 \Psi_{i,j} = \dfrac{1}{6\,h^2}$
\begin{tabular}{|c|c|c|} \hline
 1 &  4 &  1 \\ \hline
 4 &-20 &  4 \\ \hline
 1 &  4 &  1 \\ \hline
\end{tabular}
$\Psi_{i,j} + O(h^2).$
\end{tabular}
\end{equation}
However, the nine-point Laplacian of Eq.~(\ref{9ptlap}), while more accurate, is still second-order.  In fact, it can easily be shown that there \emph{cannot} exist a fourth-order nine-point Laplacian operator \cite{ME_DISS} (this fact should not be confused with the well-known fourth-order nine-point scheme for the Laplace and Poisson equations \cite{FD_PDE_BOOK}).  Even though the nine-point Laplacian of Eq.~(\ref{9ptlap}) is a more accurate second-order approximation, in order to minimize the amount of computation needed for the 2SHOC scheme, we will only utilize the standard five-point Laplacian of Eq.~(\ref{dx2full2d}) for the required second-order derivatives.  Once again, `double-differencing' is used to obtain a second-order accurate approximation for the truncation term's fourth derivatives.  However, in this case, the result yields two variations of the 2SHOC scheme.

The first is a direct parallel to the one-dimensional 2SHOC scheme in which the second derivative in the $x$ and $y$ directions are approximated with second-order central-differencing.  These are then used to form second-order approximations to the fourth derivatives in the truncation terms in Eq.~(\ref{dx2full2d}).  After simplification, this results in the 2SHOC scheme
\begin{alignat}{3}
&1)\qquad &D^x_{i,j} &= \delta_{x}^2\Psi_{i,j} = \frac{\Psi_{i+1,j} - 2\Psi_{i,j} +  \Psi_{i-1,j}}{h^2}, \label{2d2shocMs1}
\\
&\qquad &D^y_{i,j} &= \delta_{y}^2\Psi_{i,j} = \frac{\Psi_{i,j+1} - 2\Psi_{i,j} +  \Psi_{i,j-1}}{h^2}, \notag 
\\
&2)\qquad &\nabla^2\Psi_{i,j} &\approx \frac{7}{6}\left(D_{i,j}^x + D_{i,j}^y\right) - \frac{1}{12}\left(D_{i+1,j}^x +  D_{i-1,j}^x + D_{i,j+1}^y +  D_{i,j-1}^y\right).
\label{2d2shocMs2}
\end{alignat}
The 2SHOC scheme as described in Eqs.~(\ref{2d2shocMs1}) and (\ref{2d2shocMs2}) require two storage arrays ($D_x$ and $D_y$) in addition to that for $\Psi$.  In large-scale computations, memory use can be a major bottleneck and so it is useful to limit the amount of extra storage required as much as possible.  We find that the two-dimensional 2SHOC scheme can be re-formulated to only require \emph{one} extra storage matrix, at the cost of requiring more computations (a trade-off that is evaluated in Sec.~\ref{s:ccompare})

To formulate the lower-storage version of the two-dimensional 2SHOC, the standard five-point second-order finite-difference approximation to the two-dimensional Laplacian is used for the first stage (stored in $D$), and then the result is combined with a second-order cross-derivative stencil to yield the second-order approximations to the fourth derivative truncation terms in Eq.~(\ref{dx2full2d}). First, we note that
\[
\nabla^2\left(\nabla^2 \Psi\right) = \frac{\partial^2}{\partial x^2}\,\nabla^2\Psi + \frac{\partial^2}{\partial y^2}\,\nabla^2\Psi = \frac{\partial^4 \Psi}{\partial x^4} + \frac{\partial^4 \Psi}{\partial y^4} + 2 \frac{\partial^4 \Psi}{\partial x^2\,\partial y^2},
\]
in which case the fourth derivative truncation terms in Eq.~(\ref{dx2full2d}) can be written as
\begin{equation}
\label{trunc}
\frac{\partial^4 \Psi}{\partial x^4} + \frac{\partial^4 \Psi}{\partial y^4} = \frac{\partial^2}{\partial x^2}\,\nabla^2\Psi + \frac{\partial^2}{\partial y^2}\,\nabla^2\Psi - 2 \frac{\partial^4 \Psi}{\partial x^2\,\partial y^2}.
\end{equation}
The cross-derivative in Eq.~(\ref{trunc}) is known to have the nine-point second-order compact stencil \cite{HOC_PHD}
\begin{equation}
\label{cdiv}
\begin{tabular}{c} 
$\left.\dfrac{\partial^4 \Psi}{\partial x^2\,\partial y^2}\right|_{i,j} = \dfrac{1}{h^4}$
\begin{tabular}{|c|c|c|} \hline
 1 & -2 &  1 \\ \hline
-2 &  4 & -2 \\ \hline
 1 & -2 &  1 \\ \hline
\end{tabular}
$\Psi_{i,j} + O(h^2).$
\end{tabular}
\end{equation}
Part of the stencil of Eq.~(\ref{cdiv}) can be written in terms of the five-point second-order approximation to the Laplacian stored in $D$, in which case the single-storage 2SHOC scheme (after simplification) is given as
\begin{alignat}{3}
&\begin{tabular}{ll} 
$1)$ & $D_{i,j} =\delta_{x}^2\Psi_{i,j} + \delta_{y}^2\Psi_{i,j} = \dfrac{1}{h^2}$ 
\begin{tabular}{|c|c|c|} \hline
  &  1 &   \\ \hline
1 & -4 & 1 \\ \hline
  &  1 &  \\ \hline
\end{tabular}
$\Psi_{i,j}$
\end{tabular} \label{2d2shocs1}
\\
\;&\; \notag
\\
&\begin{tabular}{ll} 
$2)$ & $\nabla^2\Psi_{i,j} \approx -\dfrac{1}{12}$
\begin{tabular}{|c|c|c|} \hline
  &   1 &   \\ \hline
1 & -12 & 1 \\ \hline
  &   1 &   \\ \hline
\end{tabular}
$D_{i,j} + \dfrac{1}{6\,h^2}$  
\begin{tabular}{|c|c|c|} \hline
1 &    & 1 \\ \hline
  & -4 &   \\ \hline
1 &    & 1 \\ \hline
\end{tabular}
$\Psi_{i,j}.$
\end{tabular} \label{2d2shocs2}
\end{alignat}

As in the one-dimensional case, both formulations of the 2SHOC schemes in two dimensions are computationally equivalent to the non-compact fourth-order stencil
\begin{equation}
\label{2dnonhoc}
\begin{tabular}{c} 
$\nabla^2\Psi_{i,j} = -\dfrac{1}{12\,h^2}$
\begin{tabular}{|c|c|c|c|c|} \hline
  &     &   1 &     &   \\ \hline
  &     & -16 &     &   \\ \hline
1 & -16 &  60 & -16 & 1 \\ \hline
  &     & -16 &     &   \\ \hline
  &     &   1 &     &   \\ \hline
\end{tabular}
$\Psi_{i,j}+ O(h^4).$
\end{tabular}
\end{equation}.

\subsection{Three-Dimensional Formulation}
\label{sec:2shoc3d}
In three dimensions, the Laplacian operator applied to $\Psi$ can be represented as
\begin{alignat}{3}
\label{dx2full3d}
\nabla^2\Psi_{i,j,k} &= \left. \frac{\partial^2 \Psi}{\partial x^2}\right|_{i,j,k} + \left. \frac{\partial^2 \Psi}{\partial y^2}\right|_{i,j,k} + \left. \frac{\partial^2 \Psi}{\partial z^2}\right|_{i,j,k} 
\\
 &= \delta_x^2 \Psi_{i,j,k} + \delta_y^2 \Psi_{i,j,k} + \delta_z^2 \Psi_{i,j,k} 
\notag 
\\
&\;\; -\frac{h^2}{12}\left( \left. \frac{\partial^4 \Psi}{\partial x^4}\right|_{i,j,k} + \left. \frac{\partial^4 \Psi}{\partial y^4}\right|_{i,j,k} + \left. \frac{\partial^4 \Psi}{\partial y^4}\right|_{i,j,k}\right) + O(h^4).\notag
\end{alignat}
There are again two formulations of the 2SHOC scheme which trade off storage requirement versus number of computations.  The first formulation takes \emph{three} storage matrices ($D^x,D^y,D^z$), while, as in the two-dimensional case, the other formulation takes only \emph{one} ($D$).  The three-storage 2SHOC directly follows from the two-dimensional version and is defined as
\begin{alignat}{3}
\label{3d2shocMs1}
&1)   \qquad &D^x_{i,j,k}          &= \delta_{x}^2\Psi_{i,j,k} = \frac{\Psi_{i+1,j,k} - 2\Psi_{i,j,k} +  \Psi_{i-1,j,k}}{h^2}
\\
& \qquad &D^y_{i,j,k}          &= \delta_{y}^2\Psi_{i,j,k} = \frac{\Psi_{i,j+1,k} - 2\Psi_{i,j,k} +  \Psi_{i,j-1,k}}{h^2} \notag
\\
& \qquad &D^z_{i,j,k}          &= \delta_{z}^2\Psi_{i,j,k} = \frac{\Psi_{i,j,k+1} - 2\Psi_{i,j,k} +  \Psi_{i,j,k-1}}{h^2} \notag
\\
&2)   \qquad &\nabla^2\Psi_{i,j,k} &\approx \frac{7}{6}\left(D_{i,j,k}^x + D_{i,j,k}^y + D_{i,j,k}^z\right) 
\label{3d2shocMs2}
\\
&\;\; \qquad &\qquad               &- \frac{1}{12}\left(D_{i+1,j,k}^x +  D_{i-1,j,k}^x + D_{i,j+1,k}^y +  D_{i,j-1,k}^y +  D_{i,j,k+1}^z +  D_{i,j,k-1}^z\right). \notag
\end{alignat}

The single storage formulation requires the use of the cross-derivatives of $\nabla^2(\nabla^2\Psi)$ as before.  We write the fourth derivative truncation terms of Eq.~(\ref{dx2full3d}) as
\begin{equation}
\label{trunc3}
\frac{\partial^4 \Psi}{\partial x^4} + \frac{\partial^4 \Psi}{\partial y^4} + \frac{\partial^4 \Psi}{\partial z^4}= \nabla^2 \left(\nabla^2 \Psi\right) - 2 \left(\frac{\partial^4 \Psi}{\partial x^2\,\partial y^2} + \frac{\partial^4 \Psi}{\partial x^2\,\partial z^2} + \frac{\partial^4 \Psi}{\partial y^2\,\partial z^2}\right),
\end{equation}
therefore the cross-derivatives in Eq.~(\ref{trunc3}) can be approximated to second-order accuracy using nine-point compact stencils in each direction which contain the same coefficients as the two-dimensional stencil of Eq.~(\ref{cdiv}).  

After extracting $D$ terms out of the cross derivatives and simplifying, the single-storage 2SHOC scheme in three dimensions becomes
\begin{alignat}{5}
&\begin{tabular}{ll} 
$1)$ & $D_{i,j,k} = \dfrac{1}{h^2}\left(\;
\begin{tabular}{|c|c|c|} \hline
  &   &   \\ \hline
  & 1 &   \\ \hline
  &   &   \\ \hline
\end{tabular}
\;\Psi_{i,j+1,k} +
\begin{tabular}{|c|c|c|} \hline
  &  1 &   \\ \hline
1 & -6 & 1 \\ \hline
  &  1 &   \\ \hline
\end{tabular}
\;\Psi_{i,j,k} +
\begin{tabular}{|c|c|c|} \hline
  &   &   \\ \hline
  & 1 &   \\ \hline
  &   &   \\ \hline
\end{tabular}
\;\Psi_{i,j-1,k}\;\right),$
\end{tabular}
\label{3d2shocs}
\\
\;&\; \notag
\\
&\begin{tabular}{ll}
$2)$ & $\nabla^2\Psi_{i,j,k} \approx -\dfrac{1}{12}\left(\;
\begin{tabular}{|c|c|c|} \hline
  &   &   \\ \hline
  & 1 &  \\ \hline
  &   &   \\ \hline
\end{tabular}
\;D_{i,j+1,k} +
\begin{tabular}{|c|c|c|} \hline
  &   1 &   \\ \hline
1 & -10 & 1 \\ \hline
  &   1 &   \\ \hline
\end{tabular}
\;D_{i,j,k} +
\begin{tabular}{|c|c|c|} \hline
  &   &   \\ \hline
  & 1 &   \\ \hline
  &   &   \\ \hline
\end{tabular}
\;D_{i,j-1,k} \; \right)$
\end{tabular}
\label{3d2shocs2} \\
&\begin{tabular}{lll}
$\qquad$ & $\qquad$ & $+ \dfrac{1}{6\,h^2}\left(\;
\begin{tabular}{|c|c|c|} \hline
  & 1 &   \\ \hline
1 &   & 1 \\ \hline
  & 1 &   \\ \hline
\end{tabular} 
\;\Psi_{i,j+1,k} +
\begin{tabular}{|c|c|c|} \hline
1 &     & 1  \\ \hline
  & -12 &    \\ \hline
1 &     & 1  \\ \hline
\end{tabular}
\;\Psi_{i,j,k} +
\begin{tabular}{|c|c|c|} \hline
  & 1 &   \\ \hline
1 &   & 1 \\ \hline
  & 1 &   \\ \hline
\end{tabular}
\;\Psi_{i,j-1,k} \; \right).$
\end{tabular} \notag
\end{alignat}
Again, both formulations of the 2SHOC scheme in three dimensions are computationally equivalent to the standard non-compact fourth-order stencil.
 
\section{Implementation of the 2SHOC Schemes for solving the NLSE}
\label{s:2shocrk4}
In this section, we show the implementation of the 2SHOC schemes into explicit methods for simulating the NLSE.  The NLSE is a good model for use with the 2SHOC schemes, as it is very relevant in current research, and as mentioned in Sec.~\ref{s:intro}, its nonlinear terms make the use of explicit time-stepping desirable.  We emphasize however, that the NLSE is being used as one example of implementing the 2SHOC scheme for the Laplacian, and that the 2SHOC scheme is independent of any additional terms in the governing PDE, and can therefore be used with a multitude of models.

There exists many numerical methods for integrating the NLSE (see reviews of Refs.~\cite{FD_NLSE_REVIEW,NUM_NLSE_METHODS_1984}) including the widely-used split-step Fourier (SSF) method \cite{NLS_SPLIT_STEP}.  Although the SSF method can be very accurate and efficient, it is not used in all cases since it is typically formulated to be only second-order accurate in time, can be computationally expensive in three-dimensional settings, is somewhat difficult to parallelize (see Ref.~\cite{NUM_NLSE_SSF_PAR}), and its required use of periodic boundary conditions can make in unsuitable for certain problems (such as homogeneous background density solutions such as the dark soliton example shown in Sec.~\ref{s:num}).  Other split-step methods which use finite-difference for the linear part of the solve are also used (see for example, the HOC split-step method of Ref.~\cite{HOC_1DNLSSS}) but those scheme typically require the use of tri-diagonal solvers for each dimension to be stable which can complicate their parallelization.  Although for simplicity and time-accuracy we use a standard finite-difference approach for time-stepping the NLSE with the 2SHOC scheme, we note that since the 2SHOC scheme for the Laplacian is independent of the time-stepping methodology, it could also be implemented into a split-step method for the NLSE.

\subsection{Explicit Time Integration using Runge-Kutta}
To formulate a fully explicit scheme for simulating the NLSE, a method of lines approach is utilized, where the 2SHOC is used for the spatial Laplacian and the resulting semi-discrete system of ODEs is integrated using an explicit time-stepping scheme.  As shown in Refs.~\cite{ME_RK4STB} and \cite{RK4_2CNLSE_STB}, the first-order ($O(k)$ where $k$ is the time-step) forward difference and the second-order ($O(k^2)$) Runge-Kutta scheme (Heun's method) for simulating the NLSE are unconditional \emph{unstable} (unless, as shown in Ref.~\cite{FD_NLSEGP_EXP_STB}, the real and imaginary parts of the NLSE are computed in staggered time steps).  Therefore, for our implementation, a higher-order Runge-Kutta method is required and is chosen to be the standard fourth-order Runge-Kutta (RK4) scheme \cite{RK4}.  To facilitate our computational cost comparisons in Sec.~\ref{s:ccompare}, we write the RK4 scheme applied to the NLSE algorithmically as 
\begin{alignat}{5}
\label{RK4}
&1)\; K_{\mbox{\scriptsize tot}} = F(\Psi^n) 
&\qquad &6)\;  K_{\mbox{\scriptsize tmp}} = F(\Psi_{\mbox{\scriptsize tmp}})
\\
&2)\;  \Psi_{\mbox{\scriptsize tmp}} = \Psi^n + \frac{k}{2}\, K_{\mbox{\scriptsize tot}}   
&\qquad &7)\;  K_{\mbox{\scriptsize tot}} = K_{\mbox{\scriptsize tot}} + 2\, K_{\mbox{\scriptsize tmp}} \notag
\\
&3)\;  K_{\mbox{\scriptsize tmp}} = F(\Psi_{\mbox{\scriptsize tmp}})  
&\qquad &8)\;  \Psi_{\mbox{\scriptsize tmp}} = \Psi^n + k\,K_{\mbox{\scriptsize tmp}} \notag
\\
&4)\;  K_{\mbox{\scriptsize tot}} = K_{\mbox{\scriptsize tot}} + 2\, K_{\mbox{\scriptsize tmp}} 
&\qquad &9)\;  K_{\mbox{\scriptsize tmp}} = F(\Psi_{\mbox{\scriptsize tmp}}) \notag
\\
&5)\;  \Psi_{\mbox{\scriptsize tmp}} = \Psi^n + \frac{k}{2}\, K_{\mbox{\scriptsize tmp}} 
&\qquad &10)\;  \Psi^{n+1} = \Psi^n + \frac{k}{6}\, (K_{\mbox{\scriptsize tot}} + K_{\mbox{\scriptsize tmp}}), \notag 
\end{alignat}   
where
\[
F(\Psi) = \frac{\partial \Psi}{\partial t} = i\left[a\nabla^2\Psi + \left(s\,|\Psi|^2 - V(\bf{r})\right)\Psi\right],
\]
and $\nabla^2\Psi$ inside $F(\Psi)$ is computed with either a standard central differencing (CD) or the 2SHOC scheme.  We denote the combined time-space schemes as RK4+CD and RK4+2SHOC respectively.

The RK4 time-stepping requires three storage matrices ($K_{\mbox{\scriptsize tmp}}, K_{\mbox{\scriptsize tot}}, \Psi_{\mbox{\scriptsize tmp}}$) in addition to the storage for the solution $\Psi$ and external potential $V({\bf r})$ (if $V({\bf r})$ is chosen to be stored rather than computed at each evaluation).  Lower-storage fourth-order Runge-Kutta schemes with comparable accuracy have been developed (a $3N$ in Ref.~\cite{RK4_3N}, and a five-stage $2N$ in Ref.~\cite{RK4_2N_5STAGE} which requires an additional function evaluation, and whose coefficients are numerically derived), but for simplicity we use the classic RK4 of Eq.~(\ref{RK4}).  Using the 2SHOC scheme inside $F(\Psi)$ requires additional storage, which we discuss in Sec.~\ref{s:ccompare}.

The RK4+2SHOC scheme is conditionally stable, in that the size of the time-step is limited by a bound based on the spatial step $h$ (see Sec.~\ref{s:stb} for the stability bounds of the RK4+2SHOC scheme applied to the NLSE).  Conditional stability is one of the only drawbacks of using explicit schemes.  However, even though implicit schemes are usually unconditionally stable, error requirements and algorithm complexity often make the explicit schemes more efficient even taking the time-step limitations into account.  

\subsection{Boundary Conditions}
\label{s:bc}
The use of proper boundary conditions is very important when performing numerical simulations.  The 2SHOC scheme requires two boundary conditions, one for each step.  The first step is a boundary condition on the Laplacian of the wavefunction (or on the separate spatial derivatives if the high-storage 2SHOC is being used).  The second step either requires a boundary condition for the fourth spatial derivatives, or in many cases, the scheme for the overall model equations dictate the boundary condition (for example, in our RK4 scheme for the NLSE, the boundary condition of the second step of the 2SHOC is not required, as the boundary condition of the time derivative of the wavefunction ($\Psi_t$) overwrites any condition that would be in the 2SHOC scheme).

As mentioned in Ref.~\cite{HOC_TIME}, it appears that in general, compact boundary conditions are not possible to realize for one-stage HOC schemes.  However, even when the use of one-sided differencing (or any other non-compact technique) is necessary, it does not negate the boundary advantage of using HOC schemes (since the advantage is not having to alter the scheme \emph{near} the boundaries).  There are some boundary condition techniques that can be seen as compact, and hence fit very well into the framework of HOC schemes. These conditions are a Dirichlet ($\Psi=\mbox{const}$) and modulus-squared Dirichlet (MSD) ($|\Psi|^2 = \mbox{const}$) \cite{ME_MSD} boundary conditions.  

In many scenarios, one would ideally want a transparent boundary condition.  However, such conditions can be complicated to implement (see Ref.~\cite{MSD_NLSE_TRANS_REV} for a review of different approaches).  For many problems, an easy alternative is to make the computational grid large enough so that one can reasonably use Dirichlet conditions at the boundaries.  When $s>0$ (attractive or focusing nonlinearity) in the NLSE, most dynamics are of areas of density which decay into a surrounding background area of zero-density, in which case a Dirichlet boundary condition of zero is useful.  In order to use the Dirichlet boundary condition with the 2SHOC scheme, the boundary condition on the Laplacian must be computed in the first step.  This is easily found by setting $\Psi_t=0$ in Eq.~(\ref{NLSE}) yielding
\begin{equation}
\label{dbc_lap}
\nabla^2\Psi_b =  -\frac{(s|\Psi_b|^2 - V_b)\Psi_b}{a},
\end{equation}
where $b$ represents a boundary grid point. In the second stage of the 2SHOC, the boundary condition of Eq.~(\ref{dbc_lap}) would also be used for the overall Laplacian boundary.  In the RK4+2SHOC method for the NLSE, the second stage boundary condition would simply apply the Dirichlet condition directly as
\begin{equation}
\label{dbc_ut}
\left. \frac{\partial \Psi}{\partial t}\right|_b =  0.
\end{equation}
 
In the case where $s<0$ (repulsive or defocusing nonlinearity) and $V({\bf r})=0$ in the NLSE, most interesting dynamics are within areas of low densities which increase towards an area with a constant background density.  In such a case, one cannot use standard Dirichlet boundary conditions since the constant boundary is in the modulus-squared of the wavefunction, not at a single real and imaginary value.  To solve this problem, we recently have developed a new modulus-squared Dirichlet boundary condition that accurately simulates a constant density at the boundaries \cite{ME_MSD}.  This boundary condition is described as
\begin{equation}
\label{msdbc_ut}
\Psi_{t,b} \approx i\,\mbox{Im}\left[\frac{\Psi_{t,b-1}}{\Psi_{b-1}}\right]\,\Psi_b,
\end{equation}
where $\Psi_t$ is the temporal derivative of $\Psi$ and $b-1$ represents the closest internal point relative to a boundary grid point in the normal direction. The standard form of the MSD boundary condition requires to first compute all interior points and then compute the boundaries, making the general form of the MSD not compatible with implicit schemes (see Ref.~\cite{ME_MSD} for details on possible implicit scheme implementation strategies of the MSD).  

The MSD boundary condition can be used for any time-dependent complex PDE and in multiple dimensions.  For our application to the NLSE, we require boundary conditions for each step of the 2SHOC scheme.  For the first step, one can substitute the NLSE in the MSD boundary condition and get an MSD condition for the Laplacian given by
\begin{equation}
\label{msdbc_lap}
\nabla^2\Psi_b \approx \left[ \mbox{Im}\left(i\, \frac{\nabla^2\Psi_{b-1}}{\Psi_{b-1}}\right) + \frac{1}{a}\,\left(N_{b-1} - N_b\right)\right]\Psi_b,
\end{equation}
where
\begin{equation}
\label{nbnb1}
N_b = s\,|\Psi_b|^2 - V_b, \qquad N_{b-1} = s\,|\Psi_{b-1}|^2 - V_{b-1},
\end{equation}
while for the second step, the boundary condition of Eq.~(\ref{msdbc_ut}) can (in the RK4+2SHOC scheme) be used directly.

It should be noted that the MSD boundary condition cannot be directly applied to the multi-storage version of the 2SHOC schemes because the MSD computes the boundary value of $\nabla^2 \Psi$, not the individual second spatial derivatives needed for Eq.~(\ref{2d2shocMs1}) or Eq.~(\ref{3d2shocMs1}).  This can be overcome by rearranging the boundary condition taking advantage of the cross derivative stencils.   However, since we will be only using the single-storage 2SHOC schemes in our examples, we do not show this alteration to the MSD boundary condition here.

There are numerous other boundary conditions one might use for the NLSE including other compact-friendly ones such as Laplacian-zero ($\nabla^2 \Psi = 0$) and periodic.  However, in this paper we have limited ourselves to Dirichlet and MSD, as they will be used in the numerical tests of Sec.~\ref{s:num}.

\subsection{Stability}
\label{s:stb}
Since the most important drawback in using explicit schemes is that they are conditionally stable, it is very important to limit this deficiency as much as possible by computing the stability bound to determine the largest time-step that is usable.
 
Since the 2SHOC schemes are algebraically equivalent to the standard non-compact fourth-order schemes, the stability bounds on the time-step should not adversely be affected by the use of the compact schemes.  However, the form of the boundary condition implementation in the 2SHOC schemes will change the way the near-boundary points are computed.  In Ref.~\cite{ME_RK4STB} we used an extension to the methodology employed in Ref.~\cite{RK4_2CNLSE_STB} to conduct an analysis of the stability bounds for simulating the multidimensional NLSE with RK4 for various boundary conditions including those in Sec.~\ref{s:bc}.  The analysis was done for both the standard second-order central differencing, as well as fourth-order differencing using the 2SHOC formulation.  The results utilizing the two boundary conditions discussed in Sec.~\ref{s:bc} are summarized as follows:
    
In the linear case where $s=0$ and with no external potential ($V({\bf r})=0$), utilizing Dirichlet boundary conditions, the stability bound on the time-step $k$ using a fourth-order central-difference scheme (with interior points computed in the 2SHOC methodology) in a $d$-dimensional setting is
\begin{equation}
\label{kd2shoclin}
k_{\mbox{\scriptsize linear}} < \left(\frac{3}{4}\right) \frac{h^2}{d\,\sqrt{2}\, a},
\end{equation}
which is just $3/4$ of the stability bound in the second-order central differencing case.

In the general NLSE case, one can get a linearized stability bound by treating the nonlinearity $|\Psi|^2$ (and the MSD boundary condition terms) as a constant, yielding
\begin{equation}
\label{kdfull}
k < \frac{\sqrt{8}}{\max\{\lVert\vec B\rVert_{\infty},\Vert \forall L_i, L_i-\vec G\rVert_{\infty}\}}\,\frac{h^2}{a},
\end{equation}
where $\vec B$ are the boundary points as defined by Table~\ref{t:bc2}, $\vec L$ is defined as
\[
L_i = \frac{h^2}{a}\left(s|\Psi_i|^2 - V({\bf r}_i)\right),
\]
and $\vec G$ is a set of values defined in Table~\ref{t:sumresults}, determined by the dimension being used.
\begin{table}[htb] 
\centering 
\begin{minipage}{6in}
\caption{Values for $B_b$ in Eq.~(\ref{kdfull}).}
\begin{center}
\begin{tabular}{|l|c|c|c|} 
\hline
$\;$  & Dirichlet ($\Psi_b = \mbox{const}$)  & MSD ($|\Psi_b|^2 = \mbox{const}$) \\ \hline
$B_b$ & $0$  & $\dfrac{h^2}{i\,a} \dfrac{1}{\Psi_{b-1}} \left. \dfrac{\partial \Psi}{\partial t} \right|_{b-1} \,  $ \\ \hline
\end{tabular}
\end{center}
\label{t:bc2}
\end{minipage}
\end{table}
\begin{table}[htb] 
\centering 
\begin{minipage}{6in}
\caption{Values for $\vec G$ in Eq.~(\ref{kdfull}) for the one-, two-, and three-dimensional NLSE.}
\begin{center}
\begin{tabular}{|c|c|c|c|} \hline
  & $d=1$ & $d=2$ & $d=3$ \\ \hline
\; &  \;                    &\;                 &\;  \\
$\vec G$ &  $\dfrac{1}{12} \times \left\{64,63,46, \right.$ & 
$\dfrac{1}{12} \times \left\{128,127,126,110,109,\right.$   & 
$\dfrac{1}{12} \times \left\{192,191,190,189,174,\right.$ \\  
\; & 
$\left. \qquad 12,-3,-4 \right\}$                
& $\left.  \qquad 92,24,9,8,-6,-7,-8 \right\} $     
& $\qquad 173,172,156,155,138,36,21,$ \\	        
\; & \; & \;
& $\qquad \left. 20,6,5,4,-9,-10,-11,-12 \right\}$ \\ 
\hline
\end{tabular}
\end{center}
\label{t:sumresults}
\end{minipage}
\end{table}

As mentioned in Ref.~\cite{ME_RK4STB}, since these stability results are based on a linearized approximation to the full nonlinear problem, in practice, one must use a time-step that is lower than the bounds given.  In our experience, setting the time-step to be $80\%$ of the given bounds ensures stability in most cases (in one-dimensional simulations, the step size can be up to $90\%$ of the given bounds for most cases).  

\section{Numerical Results}
\label{s:num}
Here we show our numerical results for using the RK4+2SHOC scheme to integrate the NLSE.  Since the 2SHOC schemes are algebraically equivalent to the standard fourth-order non-compact schemes, they should exhibit the same order of accuracy as well.  However, because the order of computation is altered (which could introduce numeric cancellation or round-off errors), and the implementation of boundary conditions is different than when using wide stencils, numerical tests of the accuracy of the 2SHOC scheme is justified.  

To test the accuracy, we integrate the NLSE using non-trivial initial conditions and record the error versus the exact solutions for various spatial step-sizes.  To compute the error in the simulations, the real and imaginary parts of the wavefunction $\Psi$ are compared to the true solution at $100$ equal-spaced intervals throughout the simulation.  The averaged L2-norm of the wavefunction error ($E^n = \lVert \Psi^n-\Psi_{\rm exact}^n \rVert_2 /{\bf N}$, where $n$ is the current time step and ${\bf N}$ is the total number of grid points) is computed at each time interval.  Using the averaged L2-norm error is necessary since we are using a fixed domain and therefore the total L2-norm would be greater for smaller $h$ due to the increase in the number of grid points.  When the simulation is completed, we define the error of the real and imaginary parts of the whole simulation as the mean of the errors at each of the $100$ intervals ($E = \Sigma E^n / K$, where $K=100$ is the number of time intervals).  To compute the order of error, we define the error order between two simulations with spatial steps $h_1$ and $h_2$ as $O = \mbox{ln}(E_{h_1}) - \mbox{ln}(E_{h_2})$.  The overall order of the scheme is determined by taking the mean of the average of the real and imaginary error orders for each run.  For comparison, we run each simulation using the classic second-order central-differencing (CD) in space as well.  

We note that since the RK4 scheme is $O(k^4)$, and the stability bounds require that $k \propto h^2$, the errors in the simulations attributed to the time-stepping are negligible compared to those due to the spatial differencing, and therefore should not effect the accuracy tests.  

Each numerical test utilizes the single-storage formulation of the RK4-2SHOC schemes where applicable and is performed using our freely distributed\footnote{http://www.nlsemagic.com} MATLAB-interfaced GPU-accelerated code package NLSEmagic \cite{NLSEmagic}.  Specifically, we make use of the `full research scripts' and the double-precision GPU-accelerated 2SHOC integrators compiled from the version `v013' source packages for one, two. and three dimensions.  The precise configuration of the scripts used in computing the results given in this paper can be obtained from the `Reproducible Figure Package' available on the NLSEmagic website\footnotemark[\value{footnote}].

\subsection{One-Dimensional Test}
\label{s:num1d}
For the one-dimensional test of the RK4+2SHOC scheme, we use a known dynamical solution to the one-dimensional NLSE (with $V(\bf{r})=0$ and $s<0$) of a co-moving dark soliton given by \cite{BEC_RCbook}
\begin{equation}
\label{1dmds}
\Psi(x,t) = \sqrt{\frac{\Omega}{s}}\,\mbox{tanh}\left[\sqrt{\frac{|\Omega|}{2\,a}}\,(x-c\,t)\right]\,\mbox{exp}\left(i\,\left[\frac{c}{2\,a}\,x + \left(\Omega - \frac{c^2}{4\,a}\right)t\right]\right),
\end{equation}
where $\Omega$ is the frequency and $c$ is the velocity of the soliton.
We use the modulus-squared Dirichlet boundary condition ($|\Psi|^2 = \mbox{constant}$) of Eqs.~(\ref{msdbc_ut}) and (\ref{msdbc_lap}) and set the size of the computational domain large enough so that the wavefunction's density is machine epsilon ($\epsilon \approx 10^{-16}$) lower than the background density at the boundaries throughout the entire time of the simulation. We set $s=-1$, $a=1$, $c=0.5$, and $\Omega=-1$ in Eq.~(\ref{1dmds}).  The grid spatial-step $h$ is varied from $1/2$ to $1/32$.  We run the simulation to an end-time of $t=10$ with a time-step size of $k=0.0005$ (which is slightly less than the maximum stability bound for the smallest value of $h$ used) resulting in $20,000$ time steps.

\begin{figure}[htbp]
\centering
\begin{minipage}{6.5in}
\begin{center}
$\begin{array}{c}
\begin{array}{cc}
\includegraphics[width=3in]{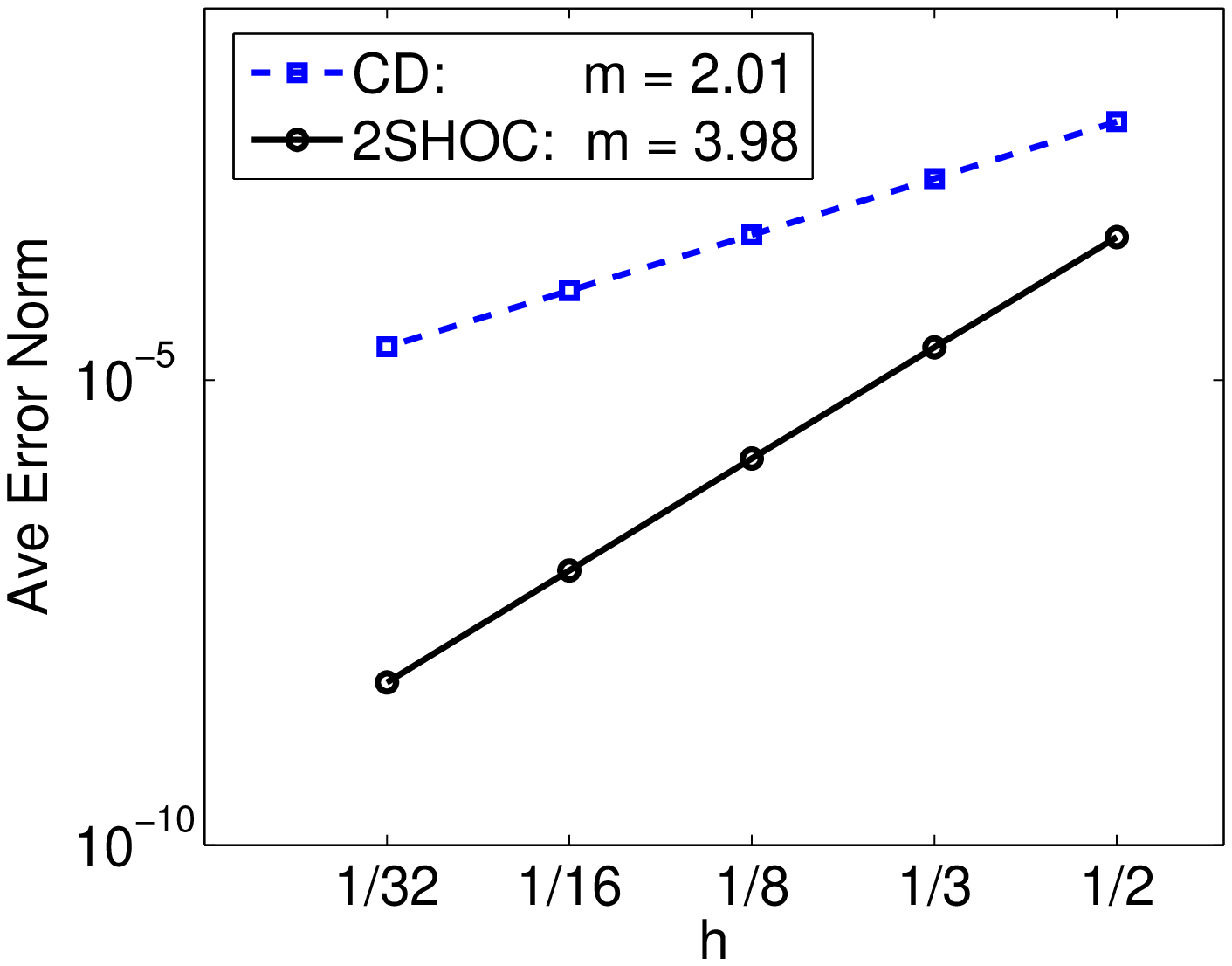} &
\includegraphics[width=3in]{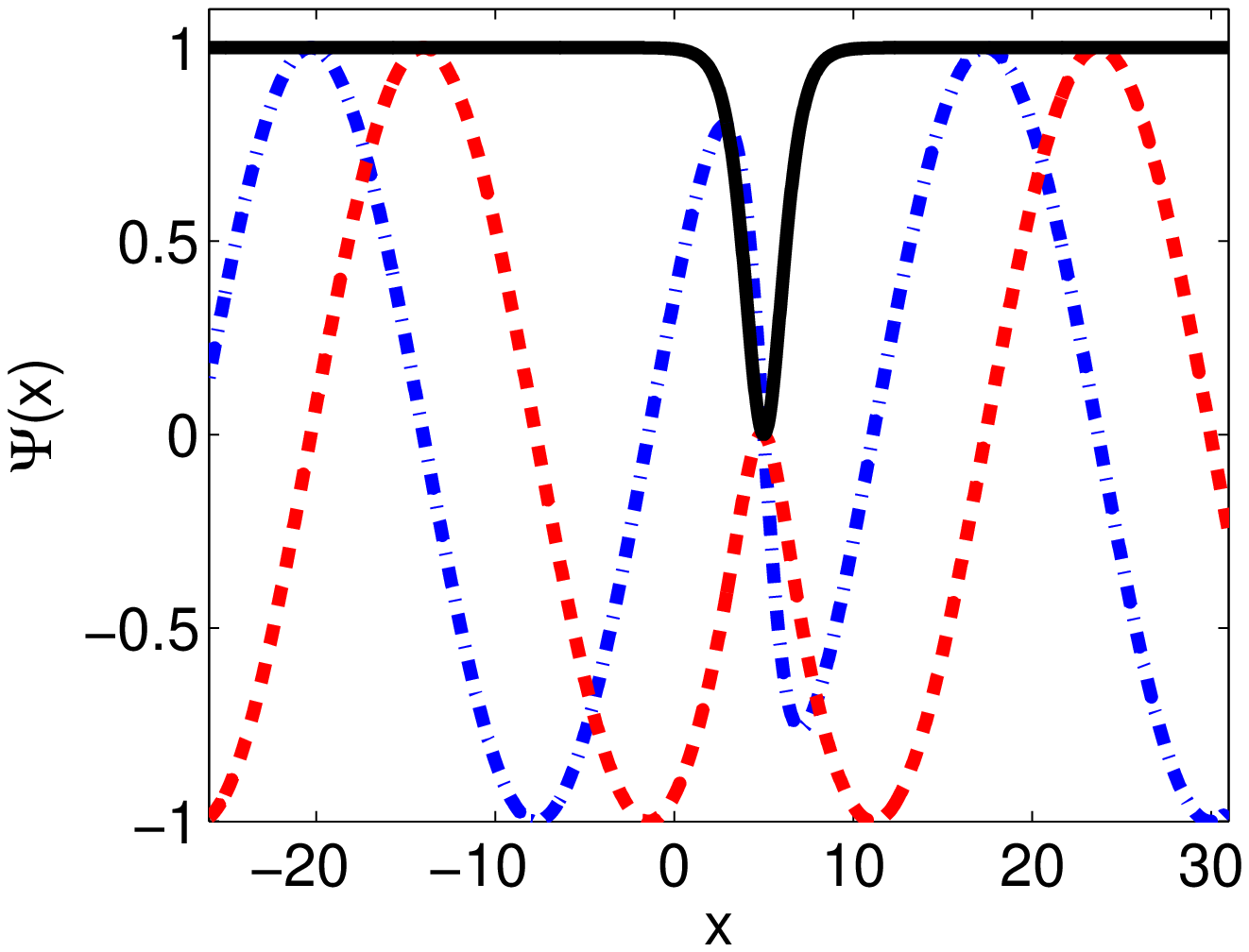}
\end{array}\\
\begin{tabular}{lccccc} \hline
 $\;$ &  $h = 1/2$  & $h = 1/4$ & $h = 1/8$ & $h = 1/16$ & $h = 1/32$  \\ \hline
\multicolumn{5}{l}{CD:}\\ \hline
Error (real) & 0.0060534 & 0.0014769 & 0.0003673 & 9.1733e-005 & 2.293e-005\\
Order (real) & -- & 2.0352 & 2.0075 & 2.0015 & 2.0002\\ \hline
Error (imag) & 0.0060159 & 0.0014675 & 0.00036496 & 9.1146e-005 & 2.2784e-005\\
Order (imag) & -- & 2.0354 & 2.0076 & 2.0015 & 2.0002\\ \hline
\multicolumn{5}{l}{2SHOC:} \\ \hline
Error (real) & 0.00034675 & 2.2662e-005 & 1.4374e-006 & 9.0208e-008 & 5.6444e-009\\
Order (real) & -- & 3.9355 & 3.9788 & 3.9941 & 3.9984\\ \hline
Error (imag) & 0.00034481 & 2.2538e-005 & 1.4297e-006 & 8.9723e-008 & 5.6141e-009\\
Order (imag) & -- & 3.9353 & 3.9786 & 3.994 & 3.9984\\ \hline
\end{tabular}
\end{array}$
\caption{(Color online) Top Left:  Overall order of accuracy ($m$) computed from simulating the one-dimensional NLSE with the initial condition of Eq.~(\ref{1dmds}) for the RK4+CD and RK4+2SHOC schemes.  Top Right:  Snap-shot of the simulation with $h=1/32$.  The solid (black) line is the modulus-squared $|\Psi|^2$, while the dot-dashed (blue) and dashed (red) lines are the real and imaginary parts of $\Psi$ respectively.  Bottom:  Table of error values and order of accuracy values for each spatial step $h$.  The parameters used for the solution are $s=-1$, $a=1$, $c=0.5$, and $\Omega=-1$.  The solution is integrated to an end-time of $t=10$ with a time-step size of $k=0.0005$.\label{f:num1d}}
\end{center}
\end{minipage}
\end{figure}

The results of the simulations are shown in Fig.~\ref{f:num1d}, where the fourth-order accuracy of the RK4+2SHOC scheme is easily observed. 

\subsection{Two-Dimensional Test}
\label{s:num2d}
\begin{figure}[htb]
\centering
\begin{minipage}{6.5in}
\begin{center}
$\begin{array}{c}
\begin{array}{cc}
\includegraphics[width=3in]{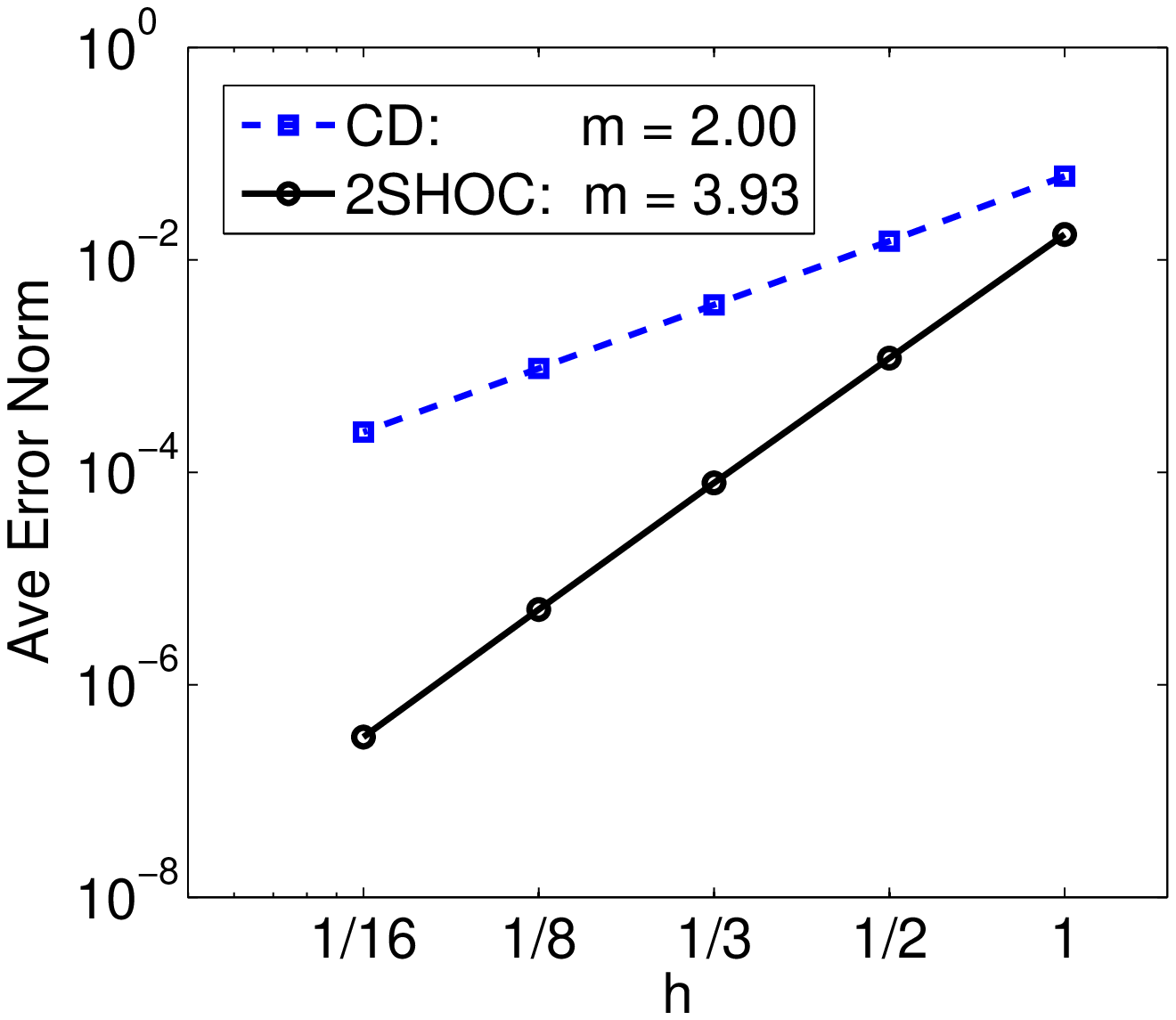} &
\includegraphics[width=3in]{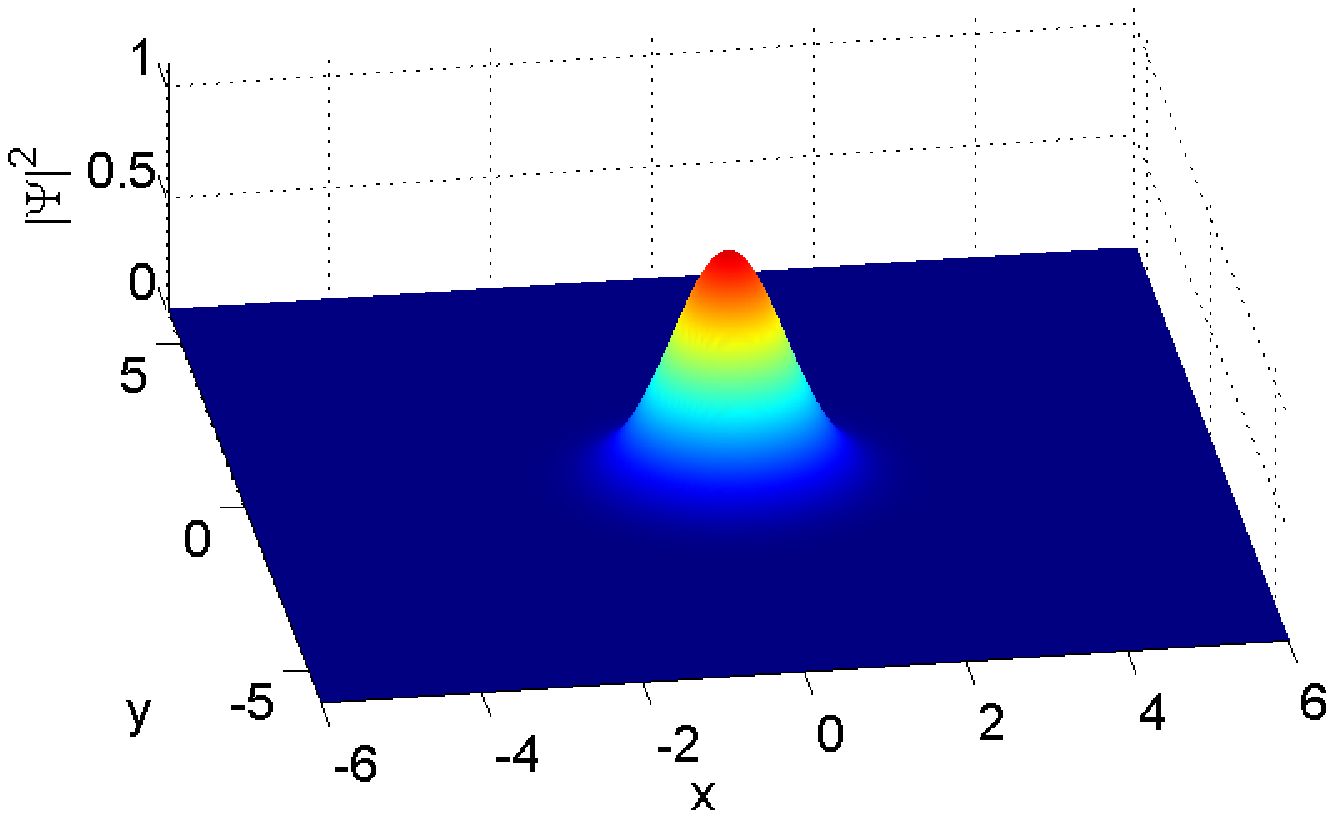}
\end{array}\\
\begin{tabular}{lccccc} \hline
 $\;$ &  $h = 1$  & $h = 1/2$ & $h = 1/4$ & $h = 1/8$ & $h = 1/16$\\ \hline
\multicolumn{5}{l}{CD:} \\ \hline
Error (real) & 0.059782 & 0.01451 & 0.0036626 & 0.00092345 & 0.00023196\\
Order (real) & -- & 2.0426 & 1.9861 & 1.9878 & 1.9932\\ \hline
Error (imag) & 0.062855 & 0.015442 & 0.0038944 & 0.00097946 & 0.00024585\\
Order (imag) & -- & 2.0252 & 1.9874 & 1.9913 & 1.9942\\ \hline
\multicolumn{5}{l}{2SHOC:} \\ \hline
Error (real) & 0.017015 & 0.0011492 & 7.7072e-005 & 4.9388e-006 & 3.1147e-007\\
Order (real) & -- & 3.8882 & 3.8982 & 3.964 & 3.987\\ \hline
Error (imag) & 0.017684 & 0.001228 & 8.2419e-005 & 5.2815e-006 & 3.3308e-007\\
Order (imag) & -- & 3.848 & 3.8972 & 3.964 & 3.987\\ \hline
\end{tabular}
\end{array}$
\caption{(Color online) Top Left:  Overall order of accuracy ($m$) computed from simulating the two-dimensional linear Schr{\"o}dinger equation with the initial condition of Eq.~(\ref{2dexpsol}) for the RK4+CD and RK4+2SHOC schemes.  Top Right:  Depiction of the modulus-squared $|\Psi|^2$ of the wavefunction with $h=1/16$ at $t=0$.  Bottom:  Table of error values and order of accuracy values for each spatial step $h$.  The parameters used for the solution are $s=0$, and $a=1$.  The solution is integrated to an end time of $t=10$ with a time-step size of $k=0.001$.\label{f:num2d}}
\end{center}
\end{minipage}
\end{figure}

There is no readily available, non-trivial, exact two-dimensional solution to the NLSE to use for order comparisons.  However, since the RK4+2SHOC scheme is explicit (and therefore does not require any special handling of the nonlinearity such as iterative methods), the accuracy of the scheme can be tested reliably in a linear setting (where $s=0$).  The chosen test problem is the Gaussian wave-packet solution
\begin{equation}
\label{2dexpsol}
\Psi(x,y,t) = \mbox{exp}\left(-\frac{x^2 + y^2}{2\,a}\right) \mbox{exp}(-i\,2\,t),
\end{equation}
where $V(x,y)$ is the external potential
\begin{equation}
\label{2dexpsolv}
V(x,y)   = \frac{x^2 + y^2}{a}.
\end{equation}
We use Dirichlet boundary conditions ($\Psi=0$) and set $a=1$.  We set the size of the computational domain large enough so that the exact solution has a value of $\Psi_b = \sqrt{\epsilon}$ (where $\epsilon \approx 10^{-16}$) at the boundaries.  The simulation time and number of intervals for error computations are the same as in the one-dimensional test in Sec.~\ref{s:num1d}.  We vary $h$ from $h=1$ to $h=1/16$ and use a time-step of $k=0.001$ resulting in $10,000$ time steps. The results are shown in Fig.~\ref{f:num2d}.
Like in the one-dimensional case, the fourth-order accuracy of the RK4+2SHOC scheme is observed.

\subsection{Three-Dimensional Test}
\label{s:num3d}
As in the two-dimensional case, there is no readily available non-trivial exact solution to test the three-dimensional NLSE.  Therefore, we again use a linear example choosing the three-dimensional analog of Eq.~(\ref{2dexpsol}) defined as
\begin{equation}
\label{3dexpsol}
\Psi(x,y,t) = \mbox{exp}\left(-\frac{x^2 + y^2 + z^2}{2\,a}\right) \mbox{exp}(-i\,3\,t),
\end{equation}
with external potential
\begin{equation}
\label{3dexpsolv}
V(x,y) = \frac{x^2 + y^2 +z^2}{a}.
\end{equation}

\begin{figure}[htb]
\centering
\begin{minipage}{6.5in}
\begin{center}
$\begin{array}{c}
\begin{array}{cc}
\includegraphics[width=3in]{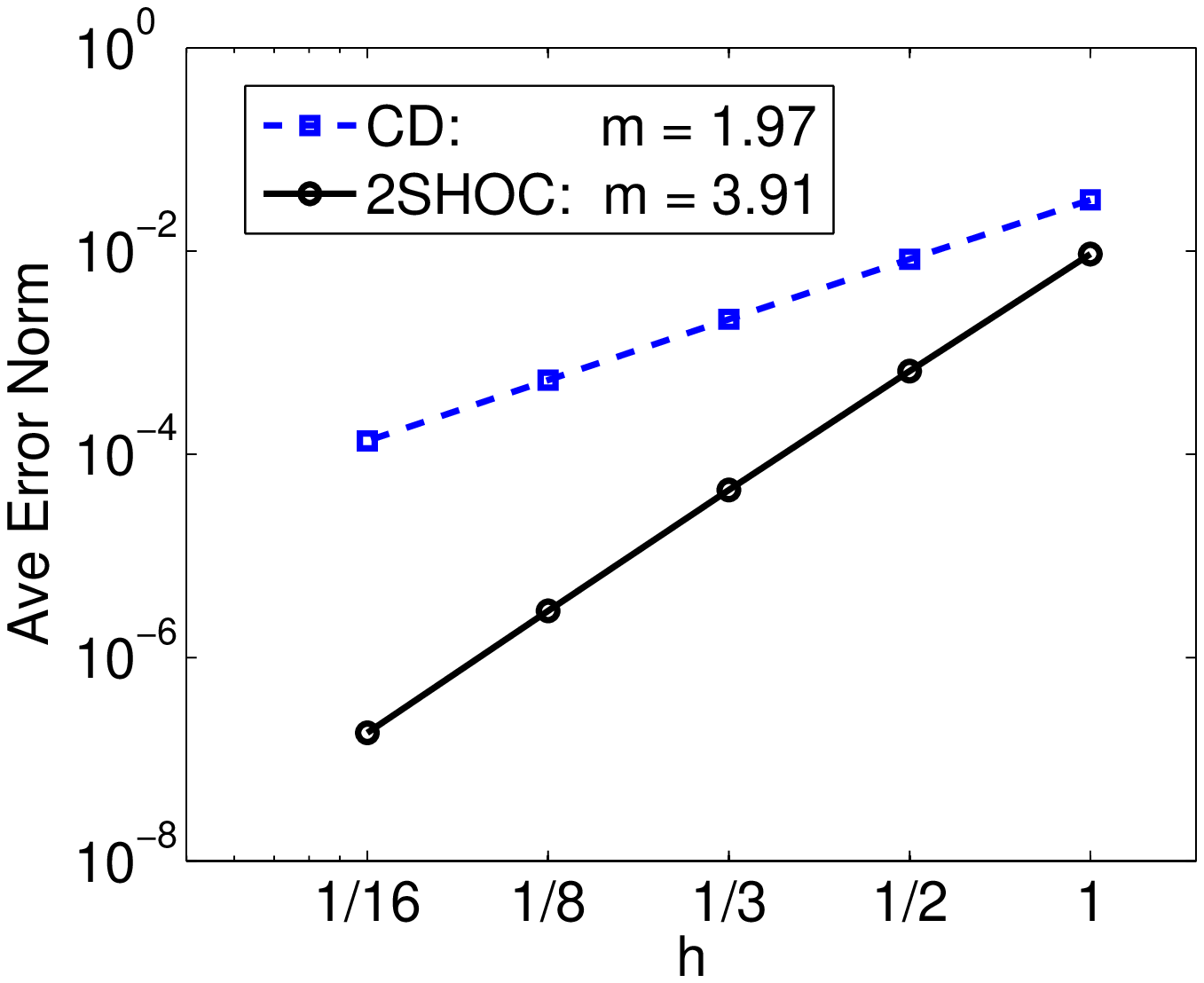} &
\includegraphics[width=3in]{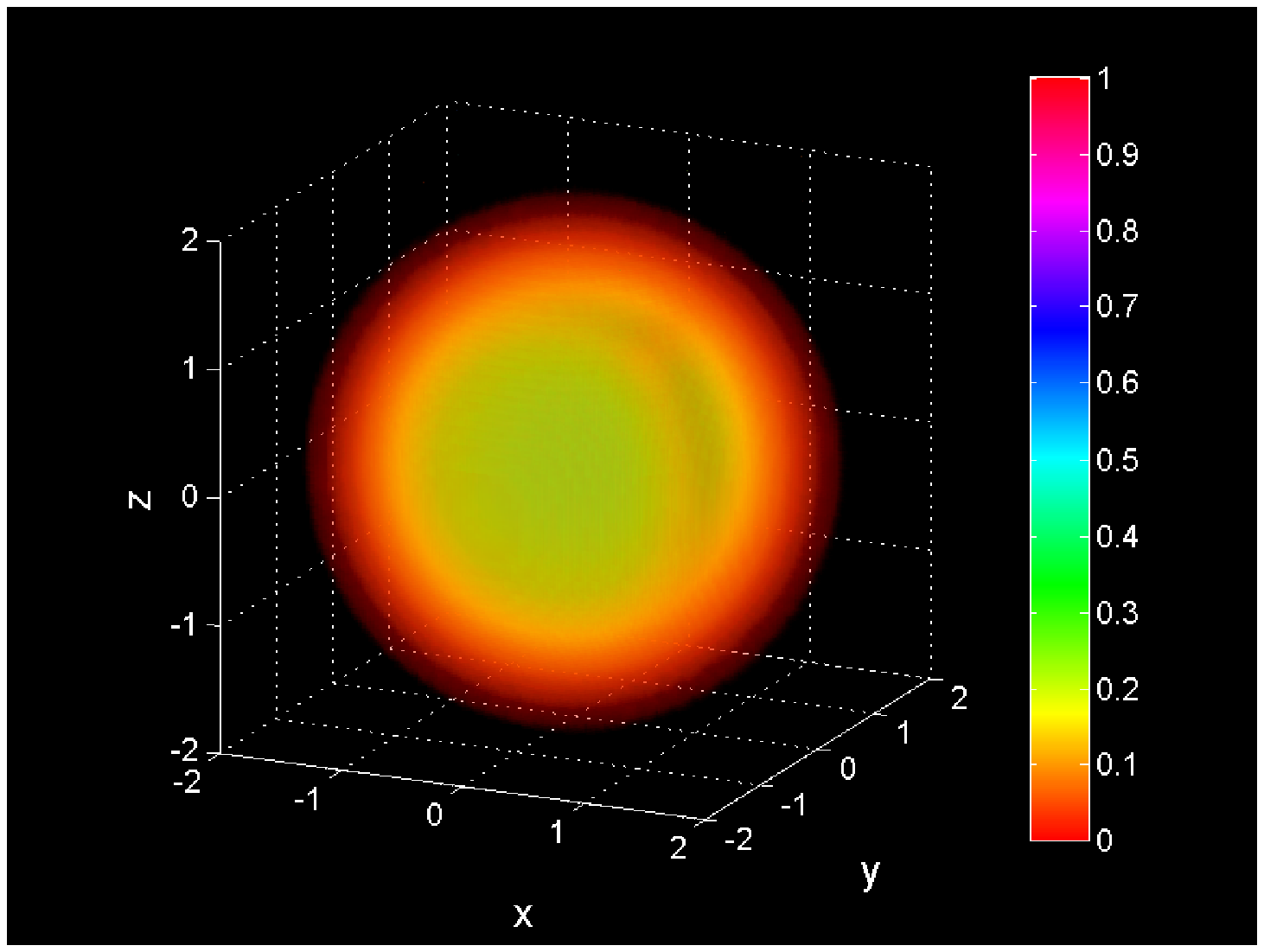}
\end{array} \\
\begin{tabular}{lccccc} \hline
 $\;$ &  $h = 1$  & $h = 1/2$ & $h = 1/4$ & $h = 1/8$ & $h = 1/16$\\ \hline
\multicolumn{5}{l}{CD:} \\ \hline
Error (real) & 0.03125 & 0.0083706 & 0.0021518 & 0.00054571 & 0.00013743\\
Order (real) & -- & 1.9004 & 1.9598 & 1.9793 & 1.9894\\ \hline
Error (imag) & 0.032541 & 0.0083009 & 0.0021026 & 0.00053144 & 0.00013374\\
Order (imag) & -- & 1.9709 & 1.9811 & 1.9842 & 1.9905\\ \hline
\multicolumn{5}{l}{2SHOC:} \\ \hline
Error (real) & 0.0093495 & 0.000666 & 4.5165e-005 & 2.9087e-006 & 1.839e-007\\
Order (real) & -- & 3.8113 & 3.8822 & 3.9567 & 3.9834\\ \hline
Error (imag) & 0.0093489 & 0.00064941 & 4.4012e-005 & 2.8344e-006 & 1.792e-007\\
Order (imag) & -- & 3.8476 & 3.8832 & 3.9568 & 3.9834\\ \hline
\end{tabular}
\end{array}$
\caption{(Color online) Top Left:  Overall order of accuracy ($m$) computed from simulating the three-dimensional linear Schr{\"o}dinger equation with the initial condition of Eq.~(\ref{3dexpsol}) for the RK4+CD and RK4+2SHOC schemes.  Top Right:  Volumetric rendering of the  modulus-squared $|\Psi|^2$ of the wavefunction with $h=1/16$ at $t=0$ (displayed on a smaller grid than used in the simulations).  Bottom:  Table of error values and order of accuracy values for each spatial step $h$.  The parameters used for the solution are $s=0$, and $a=1$.  The solution is integrated to an end time of $t=10$ with a time-step size of $k=0.0005$.\label{f:num3d}}
\end{center}
\end{minipage}
\end{figure}

As before, Dirichlet boundary conditions ($\Psi=0$) are used and $a$ is set to $1$.  The size of the computational domain is once again set to be large enough so that the exact solution has a value of $\Psi_b = \sqrt{\epsilon}$ at the boundaries.  The simulation time and number of intervals for error computations are the same as in the one-dimensional test in Sec.~\ref{s:num1d}.  The step-size is varied from $h=1$ to $h=1/16$ and a time-step of $k=0.0005$ is used resulting in $20,000$ time steps. The results are shown in Fig.~\ref{f:num3d}, where once again the fourth-order accuracy of the scheme is clearly seen.
The average order of accuracy given in Fig.~\ref{f:num3d} is $3.91$ which is slightly smaller than expected for the fourth-order scheme.  However as seen in the corresponding table, the order of accuracy starts at around $3.8$ for the first spatial step-size reduction, while in the smaller reductions, the order increases to around $3.98$.  This means that the slightly lower-order values are most likely due to the inability of the coarse grid to adequately resolve the solution, and not on the scheme itself.  Lower values of $h$ were not used due to computational memory constraints.

From the table in Fig.~\ref{f:num3d}, it is easy to illustrate the advantage of using high-order schemes for large three-dimensional problems mentioned in the beginning of this chapter --- that of being able to use a much smaller grid while maintaining the desired accuracy.  For example, in Fig.~\ref{f:num3d}, using the second-order scheme for $h=1/8$ (making the grid resolution $97\times 97 \times 97 = 912,673$ grid points) yielded an error of around $0.0005$.  Roughly the same error ($0.0006$) was found using the fourth-order scheme with a spatial-step size of $h=1/2$, requiring the grid size of only  $25\times 25\times 25 = 15,625$ grid points.  This is a $98.3\%$ reduction in the number of grid points required!  When combined with the ease of parallelization that the 2SHOC compact scheme provides, its usefulness for large three-dimensional problems is apparent.  

\section{Computation and Storage Comparisons}
\label{s:ccompare}
In this section we compare the storage and computational requirements of the 2SHOC schemes compared to the standard fourth-order non-compact equivalent schemes for one, two and three dimensions.  We count the number of operations (in terms of the number of elements in the domain, $N$) needed for each method ignoring the boundaries.  For simplicity, we also ignore the added operations needed due to $\Psi$ being complex, and treat all operations as acting on real variables.  Also, since floating-point division operations are far more computationally expensive than additions and multiplications, we record the operations required in an optimized form of the schemes, in which all divisions are only computed once (and hence ignored) by pre-computing the constant terms.  We record the number of computations and storage space required for the Laplacian operator alone using the 2SHOC, as well as when implemented into the NLSE simulations using the RK4+2SHOC schemes. 

The one-dimensional analysis is shown in Table~\ref{t:1d2shocA}. 
\begin{table}[htb] 
\centering 
\begin{minipage}{6in}
\caption{One-dimensional storage and computational cost analysis for the 2SHOC scheme compared to the equivalent non-compact scheme. The storage is given in terms of $N$, the total number of grid points.}
\begin{center}
\begin{tabular}{lcccc} \hline
 Method &  Operations  & Storage & Op Ratio & Storage Ratio \\ \hline
\multicolumn{4}{l}{Laplacian:} \\ \hline
Non-Compact & 7  & $N$ &  -- & -- \\
2SHOC       & 8  & $2N$ & 1.14 & 2 \\ \hline
\multicolumn{4}{l}{NLSE RK4 Step:} \\ \hline
Non-Compact & 4\,(7+7)+13 = 69 & $5N$ &  -- & -- \\
2SHOC       & 4\,(8+7)+13 = 73 & $6N$ & 1.06 & 1.2 \\ \hline
\end{tabular}
\end{center}
\label{t:1d2shocA}
\end{minipage}
\end{table}
We see that for the RK4+2SHOC NLSE implementation, the scheme only requires about $6\%$ more computations and $20\%$ more storage than the equivalent non-compact scheme.  This small increase in storage and computations is minor compared to the advantages of having a compact scheme.

The two-dimensional analysis is given in Table~\ref{t:2d2shocA}.  Both the double-storage version of the 2SHOC of Eqs.~(\ref{2d2shocMs1}) and (\ref{2d2shocMs2}) and the single-storage version of Eqs.~(\ref{2d2shocs1}) and (\ref{2d2shocs2}) are analyzed.  
\begin{table}[htb] 
\centering 
\begin{minipage}{6in}
\caption{Two-dimensional storage and computational cost analysis for 2SHOC schemes compared to equivalent non-compact schemes. The storage is given in terms of $N$, the total number of grid points.}
\begin{center}
\begin{tabular}{lcccc} \hline
 Method &  Operations  & Storage & Op Ratio & Storage Ratio \\ \hline
\multicolumn{4}{l}{Laplacian:} \\ \hline
Non-Compact       & 9  & $N$  & -- & -- \\
2SHOC (2X Storage)& 15 & $3N$ & 1.66 & 3 \\
2SHOC (1X Storage)& 19 & $2N$ & 2.11 & 2 \\ \hline
\multicolumn{4}{l}{NLSE RK4 Step:} \\ \hline
Non-Compact       & 4\,(9+7) +13 = 77  & $5N$ &  --   & --    \\
2SHOC (2X Storage)& 4\,(15+7)+13 = 101 & $7N$ & 1.31 & 1.4  \\
2SHOC (1X Storage)& 4\,(19+7)+13 = 117 & $6N$ & 1.52 & 1.2 \\ \hline
\end{tabular}
\end{center}
\label{t:2d2shocA}
\end{minipage}
\end{table}
Here, using the 2SHOC scheme is a greater increase in additional computations, taking 50\% more (for the 1X storage version) than the non-compact scheme in the RK4 step of the NLSE.  The two-storage version of 2SHOC reduces this to about $30\%$, but with a $20\%$ increase in the storage required when compared to the single-storage 2SHOC (which like the one-dimensional case is only $20\%$ more than the non-compact scheme).

The three-dimensional analysis is given in Table~\ref{t:3d2shocA}.
\begin{table}[htb] 
\centering 
\begin{minipage}{6in}
\caption{Three-dimensional storage and computational cost analysis for 2SHOC schemes compared to equivalent non-compact schemes. The storage is given in terms of $N$, the total number of grid points.}
\begin{center}
\begin{tabular}{lcccc} \hline
 Method &  Operations  & Storage & Op Ratio & Storage Ratio \\ \hline
\multicolumn{4}{l}{Laplacian:} \\ \hline
Non-Compact       & 14 & $N$  & --    & -- \\
2SHOC (3X Storage)& 22 & $4N$ & 1.57 & 4 \\
2SHOC (1X Storage)& 31 & $2N$ & 2.21 & 2 \\ \hline
\multicolumn{4}{l}{NLSE RK4 Step:} \\ \hline
Non-Compact       & 4\,(14+7)+13 =  97 & $5N$ &  --   & --    \\
2SHOC (3X Storage)& 4\,(22+7)+13 = 129 & $8N$ & 1.33 & 1.6 \\
2SHOC (1X Storage)& 4\,(31+7)+13 = 165 & $6N$ & 1.70 & 1.2 \\ \hline
\end{tabular}
\end{center}
\label{t:3d2shocA}
\end{minipage}
\end{table}
Here we see that the single-storage 2SHOC scheme in the RK4 step for the NLSE takes $70\%$ more computations than the non-compact scheme, and for the Laplacian operator alone, takes a little more than twice the computations required for the non-compact scheme.  The 1X-storage 2SHOC scheme for the Laplacian only takes $57\%$ more computations but also uses four times the storage.  Since the main advantage of using a compact scheme is that it is easy to implement in a parallel environment and to realize near the boundaries, the increase in computation is easily offset by the parallelism (as most parallel implementations have speed-ups of well over two).  In addition, often times, the number of computations an algorithm takes is not nearly as important as the memory bandwidth used.  Therefore, using less memory but more computations may actually perform better than using less computations and more memory.  Also, since memory capacity can be limited in certain computational environments such as GPUs, using the single-storage 2SHOC schemes allow for larger problem sizes (it is for this reason that the single-storage 2SHOC schemes were chosen for implementation in the NLSEmagic GPU code package \cite{NLSEmagic} mentioned in the beginning of this section).

\section{Conclusion}
\label{s:con}
We have formulated two-step high-order compact (2SHOC) finite-difference schemes for the Laplacian operator in one, two, and three dimensions.  Two forms of the multi-dimensional 2SHOC schemes were formulated, one requiring less storage and more computation than the other. When implementing the schemes in partial differential equation models, the schemes do not rely on other parts of the governing equation, and are therefore more general than other HOC schemes.  

The schemes are well-suited for use in explicit finite-difference schemes for time-dependent problems.  As an example, we have shown the implementation of the 2SHOC schemes into explicit schemes for solving the nonlinear Schr{\"o}dinger equation and, through numerical simulations, have shown that the schemes display the desired accuracy.  

A computation and storage analysis revealed that the 2SHOC schemes take more storage and computations (never more than roughly double) than their non-compact equivalent schemes, but that these increases are well offset due to the advantages of using a compact scheme, most notably, the easing of parallel implementations.  

\section*{Acknowledgments}
This research was supported by NSF-DMS-0806762 and the Computational Science Research Center at San Diego State University.

\def\myitemsep{5pt}
\bibliographystyle{mio}
\bibliography{2SHOC4}  
\end{document}